\documentclass[11pt]{article}
\usepackage[dvips]{graphicx}
\usepackage{a4wide}
\usepackage{amsmath}
\usepackage{amsthm}
\usepackage{amsfonts}
\usepackage{bm}
\usepackage{cases}
\usepackage{mathrsfs}
\usepackage{authblk}
\usepackage{here}

{\theoremstyle{plain}%
  \newtheorem{thm}{\bf Theorem}[section]%
  \newtheorem{lem}[thm]{\bf Lemma}%
  \newtheorem{prop}[thm]{\bf Proposition}%
  \newtheorem{claim}[thm]{\bf Claim}%
}
{\theoremstyle{remark}
  
}

\newcommand{\dC}{{\mathbb{C}}}%
\newcommand{\dN}{{\mathbb{N}}}%
\newcommand{\dR}{{\mathbb{R}}}%

\title{Periodicity of  Grover walks on generalized Bethe trees}

\author[1]{Sho KUBOTA}
\author[1]{Etsuo SEGAWA}
\author[2]{Tetsuji TANIGUCHI}
\author[1]{Yusuke YOSHIE\footnote{Corresponding author. E-mail: yusuke.yoshie.r1@dc.tohoku.ac.jp}}
\affil[1]{Graduate School of Information Sciences, Tohoku University, \newline Aoba, Sendai 980-8579, Japan}
\affil[2]{Department of Electronics and Computer Engineering, Hiroshima Institute of Technology, \newline Hiroshima 731-5193, Japan }
\date{}
\begin {document}
\maketitle{}
\begin{abstract}
We focus on the periodicity of the Grover walk on the generalized Bethe tree, which is a rooted tree such that  in each level the vertices  have the same degree. Since the Grover walk is induced by the underlying graph, its properties depend on the graph. In this paper, we say that {\it the graph induces periodic Grover walks} if and only if there exists $k \in \dN$ such that the $k$-th power of the time evolution operator becomes the identity operator. Our aim is to characterize such graphs. We give the perfect characterizations of the generalized Bethe trees which induce periodic Grover walks. 
\end{abstract}

\section{Introduction}
\subsection{Introduction and related works}
Quantum walks are introduced as  quantum versions of random walks \cite{Gudder}. The quantum walk is induced by a based graph, and the motion of the walker can be regarded as a discrete analogue of scattering of plane wave on the graph \cite{Emms06}, \cite{HSegawa17}. The state of the walker at each time is represented by an $\ell^{2}$-function on an induced Hilbert space $\mathcal{H}$, and its evolution is given by a unitary operator on $\mathcal{H}$. Reviews on quantum walks from the viewpoints of several research fields can be seen in the following books, e.g. \cite{Konno14}, \cite{Wang}, \cite{Portugal}. Studies of quantum walks have been developed for last decade, and they have been applied to various fields \cite{Chisa11}, \cite{Stra07}. In particular, the quantum search algorithm is one of the most attractive applications of quantum walks. They enable us to find marked elements on a graph faster than classical random walks \cite{Kempe}, \cite{Sze04}. Also, a perfect state transfer from the initial point to the target point can be regarded as a quantum search \cite{Godsil}, \cite{Kendon} \cite{Stefanak}. Here, we treat the periodicity of a quantum walk. The periodicity appears if there exists an integer $k$ such that the $k$-th power of the time evolution operator becomes the identity. It also implies that arbitary state at time $k$ returns the initial state. So the periodicity can be viewed as a perfect state transfer on the same point and our ultimate purpose is to characterize the graphs in which such a kind of perfect state transfer occurs.  If an underlying graph induces periodic quantum walks, then the behavior of the walker is periodic with some periods, and the sequence of the distribution is also periodic. By a spectral method, we consider the periodicity. In this paper, we focus on the periodicity of the Grover walk. The Grover walk is a kind of quantum walk, which is uniquely determined by the underlying graph and it is also related to several study fields, e.g. not only quantum searches but also an analysis of the zeta function \cite{Sato11}, an isomorphic problem for two cospectral strongly regular graphs \cite{HSegawa13}. Indeed, the spectrum of the evolution operator of the Grover walk derives from that of the isotropic random walk on the underlying graph. So the spectral analysis of the random walk can be applied to deal the periodicity of the Grover walk. 

In this paper, we say that {\it the graph induces a k-periodic Grover walk} if and only if the $k$-th power of the time evolution operator becomes the identity first ($k \in \dN$). In \cite{Higuchi13}, characterizations of some fixed finite graphs to induce a periodic Grover walk are introduced, and the results are as follows: 
\begin{itemize}
\item (Complete graphs) The Grover walk on a complete graph $K_{n}$ is periodic if and only if $n=2$, or $3$, whose periods are $2, 3$, respectively.
\item (Complete bipartite graphs) The Grover walk on a complete bipartite graph $K_{r, s}$ is periodic for any $r, s \in \dN$ with $r+s \ge 3$, whose period is $4$. 
\item (Strongly regular graphs) The Grover walk on a strongly regular graph $\mathrm{SRG}(n, k, \lambda, \mu)$ is periodic if and only if $(n, k, \lambda, \mu)=(5, 2, 0, 1), (2k,k,0,k)$, or $(3\lambda ,2\lambda ,\lambda ,2\lambda)$, whose periods are $5, 4, 12$, respectively. These graphs are nothing but $C_{5}, K_{k, k}$, and $K_{\lambda, \lambda, \lambda}$, respectively. 
\end{itemize}
Characterizations of graphs which induce  $k$-periodic Grover walks for a fixed integer $k$ are treated in \cite{Yoshie}. The results are as follows: Let $\mathscr{W}_{k}$ be the family of classes of finite simple connected graphs which induce $k$-periodic Grover walks. Then, 
\begin{itemize}
\item $\mathscr{W}_{2}=\{ P_{2} \}$,
\item $\mathscr{W}_{3}=\{ C_{3} \}$,
\item $\mathscr{W}_{4}=\{ K_{r, s} | r, s \in \dN \}$,
\item $\mathscr{W}_{5}=\{ C_{5} \}$,
\item $\mathscr{W}_{2r-1} \subset \{ \text{Odd unicycle graphs} \}$
\end{itemize}
for $r \in \dN$. For random walks, such phenomenon does not occur as long as the underlying graph is connected and non-bipartite in general. So we can say that the periodicity is an special feature of the quantum walk. All the graphs in the above lists have a special partition called {\it an equitable partition} \cite{RGodsil}. In this paper, we treat the generalized Bethe tree, which is a typical graph with an equitable partition.  

\subsection{Main result}

In this paper, we treat a special class of the tree graphs called the generalized Bethe tree. It is a rooted tree such that  in each level the vertices have the same degree. Researches on generalized Bethe trees are treated in \cite{Rojo05}, \cite{Rojo07}. Analysises of the eigenvalues of the adjacency matrix and the Laplacian matrix of generalized Bethe trees are discussed in both of the two references. In this paper, we provide an analysis of the eigenvalues of the matrix corresponding to the random walk on generalized Bethe trees and give the perfect characterizations of these graphs to induce a periodic Grover walk. Also, we give notations of some graphs. For $k \in \dN$, the $k$-subdivision graph of $G$ is denoted by $S_{k}(G)$, which is a graph obtained by putting $k-1$ vertices to each edges of $G$. Also, let $P_{l}, ST_{l}$ be the path graph and the star graph with $l$ vertices, which are special classes of the generalized Bethe tree. Then, the following Theorem is our main result in this paper.
\begin{thm}
The generalized Bethe trees $B(d(0), d(1), \cdots, d(n-1))$ which induce a periodic Grover walk are only $P_{n+1}$, $S_{k}(ST_{l})$, $S_{k}(B(1, 2, 3))$, and $S_{k}(B(s, 3))$ for $s \in \dN_{\ge 2}$, and $l, k \in \dN$, whose periods are $2n$, $4k$, $12k$, and $12k$, respectively.  \label{thm of bethe}
\end{thm}
The graphs $B(1, 2, 3)$, $B(s, 3)$ are provided in Figures \ref{fig bethe3} and  \ref{fig bethe4}. (See Section 2 for detailed definition).

This paper is organized as follows: In Section 2, we give detailed definition of the Grover walk and the generalized Bethe trees as the preliminaries. Also, we introduce some key tools to consider the periodicity of the Grover walk on the generalized Bethe trees. In Section 3, we prove our main result, Theorem \ref{thm of bethe} using these tools. We summarize our results and make a discussion in Section 4. 
\begin{figure}[h]
\hspace{2cm}
\begin{tabular}{ccc}
		\begin{minipage}[h]{0.25\linewidth}
		\centering
		\vspace{1.5cm}
		\scalebox{1.0}{
		\includegraphics[scale=1.0]{bethe-figfig.1}
		}
		\vspace{0.1cm}
              \centering
             \vspace{0.5cm}
		\caption{$B(1, 2, 3)$}
		\label{fig bethe3}
		\end{minipage}
		\hspace{1cm}
		\begin{minipage}[h]{0.25\linewidth}
		\centering
		\scalebox{1.0}{
		\includegraphics[scale=1.0]{bethe-figfig.2}
		}
		\caption{$B(s, 3)$}
		
		\label{fig bethe4}
		\end{minipage}	
\end{tabular}
\end{figure}

\section{Preliminaries}
\subsection{Definition of Grover walk}
Here, we give the definition of the Grover walk. Let $G=(V, E)$ be a finite simple connected graph with the vertex set $V$, and the edge set $E$. Then, it can be considered that each edges in $E$ have two orientations. For $uv \in E$, an arc from $u$ to $v$ is denoted by $e=(u, v)$. The origin and terminus vertex of $e$ are denoted by $o(e), t(e)$, respectively. Also, the inverse arc of $e$ is denoted by $e^{-1}$. We define $D(G)=\{ (u, v), (v, u) | uv \in E \}$, which is a set of symmetric arcs of $G$ and give a unitary matrix $U=U(G)$ indexed by $D(G)$ as follows:
\begin{equation*}
U_{e, f}=
	\begin{cases}
	2/\mathrm{deg}(t(f)) & \text{if $t(f)=o(e), e\neq f^{-1}$,}\\
        2/\mathrm{deg}(t(f))-1 & \text{if $e=f^{-1}$,}\\
        0 & \text{otherwise.}
	\end{cases}
\end{equation*}  
The time evolution operator of the Grover walk on the graph $G$ is determined by the above $U$ which is called the Grover transfer matrix. Thus, the Grover walk is induced by the underlying graph. Let $\varphi_{t} \in \ell^{2}(D(G))$ be a quantum state at time $t$. Then, the quantum state at time $t+1$, $\varphi_{t+1}$, is given by $\varphi_{t+1}=U\varphi_{t}$. A graph $G$ induces a $k$-periodic Grover walk if and only if $U^{k}=I_{|D(G)|}$ and $U^{j} \ne I_{|D(G)|}$ for every $j$ with $0<j<k$. So it immediately follows that if $G$ induces a $k$-periodic Grover walk then $\varphi_{k}$ returns to $\varphi_{0}$. For a square matrix $A$, let $\sigma(A)$ be the set of the eigenvalues of $A$. Here, one can easily show the following useful Proposition for this paper: 
\begin{prop}
A graph $G$ induces a $k$-periodic Grover walk if and only if for $U=U(G)$, $\lambda^{k}_{U}=1$ for every $\lambda_{U} \in \sigma(U)$, and there exists $\lambda_{U} \in \sigma(U)$ such that $\lambda^{j}_{U} \ne 1$ for every $j$ with $0 < j < k$.\label{prop of p}
\end{prop}
In this paper, we use spectral method in order to analyze the periodicity of the Grover walk with the above Proposition. Let $T$ be the transition matrix of $G$, that is, for $u, v \in V$, 
\begin{equation*}
(T)_{u, v}
=
\begin{cases}
	1/\mathrm{deg}(u) & \text{if $u \sim v$,}\\
	0 & \text{otherwise.}
	\end{cases}
\end{equation*}
There are results for the eigenvalues of $U$, e.g. \cite{Emms06}, \cite{HSegawa13}. In particular, we introduce the following.
\begin{lem}
(Higuchi, Segawa \cite{HSegawa17})\\
The spectrum of the Grover transfer matrix $U$ is decomposed by
\[\sigma(U)= \{ e^{\pm i \arccos{ \left( \sigma(T) \right)} } \}  \cup {\{ 1 \}}^{b_{1}}  \cup  {\{ -1 \} }^{b_{1}+1-{\bf 1_{B}}},  \]
where $b_{1}$ is the first Betti number of $G$, that is, $|E|-|V|+1$, and ${\bf 1_{B}}= 1$ if G is bipartite,  ${\bf 1_{B}}= 0$ otherwise.
\end{lem}
Therefore, if a graph $G$ induces a $k$-periodic Grover walk, then the eigenvalues of the transition matrix of $G$ should be the real part of the $k$-th root of unity. Actually, the following Zhukovskij transformation is useful to solve the problem.  
\begin{lem}
(Higuchi, Konno, Sato, Segawa \cite{Higuchi13})\\
Let $f(\lambda)$ be a monic polynomial of degree $i$ for $\lambda \in \dR$. Then, the roots of $f(\lambda)$ are the real part of the roots of unity if and only if for $z \in \dC$ with $|z|=1$, the polynomial  $(2z)^{i}f\left( (z+z^{-1})/2 \right)$ is a product of some cyclotomic polynomials. \label{lem of cyclo}
\end{lem}

\subsection{Definition of generalized Bethe trees}
Throughout this paper, we denote $[i, j]=\{ i, i+1, \cdots, j-1, j \}$ for $i, j \in \dN$ with $i < j$. We focus on a generalized Bethe tree $G$ which is a rooted tree such that in each level the vertices have the same degree. We agree that the root vertex is at level $0$, and this tree has $n+1$ levels. Let $C_{i}$ be the set of vertices of level $i$ for $i \in [0, n]$. Then, for every $i \in [0, n-1]$ the value $|N(v) \cap C_{i+1}|$ is constant whenever $v$ is in $C_{i}$, where $N(v)$ is the set of the neighbors of $v$. Then, $\{ C_{0}, C_{1}, \cdots, C_{n} \}$ is a partition of $V(G)$ and such a partition is called {\it an equitable partition} \cite{RGodsil}. So the generalized Bethe trees are tree graphs with an equitable partition. We put $d(i)=|N(v)\cap C_{i+1}|$ for $i \in [0, n-1]$, $v \in C_{i}$,  and $d(n)=0$. In this paper, a generalized Bethe tree is denoted by $B(d(0), d(1), \cdots, d(n-1))$. We provide examples of generalized Bethe trees in Figures \ref{fig bethe1} and \ref{fig bethe2}. For $v \in C_{i}$ with $i \ne 0$ and $i \ne n$, the child and parent of $v$ are vertices in $|N(v) \cap C_{i+1}|, |N(v) \cap C_{i-1}|$, respectively. 
\begin{figure}[htbp]
\hspace{1cm}
\begin{tabular}{ccc}
		\begin{minipage}[h]{0.25\linewidth}
		\centering
		\scalebox{1.0}{
		\includegraphics[scale=1.0]{bethe-fig.1}
		}
		\vspace{0.1cm}
              \centering
		\caption{$B(2, 3, 1)$}
		\label{fig bethe1}
		\end{minipage}
		\hspace{3cm}
		\begin{minipage}[h]{0.25\linewidth}
		\centering
		\scalebox{1.0}{
		\includegraphics[scale=1.0]{bethe-fig.2}
		}
		\caption{$B(5, 2)$}
		\label{fig bethe2}
		\end{minipage}		
\end{tabular}
\end{figure}
Also, we define  
\[ D_{i}=\frac{d(i)}{( d(i)+1 ) ( d(i+1)+1 )}, \]
for $i \in [1, n-1]$, and $D_{0}=1/(d(1)+1)$.  Indeed, the generalized Bethe tree can be projected on a path graph in which the levels correspond to the vertices of the path. Then $D_{i}$ is nothing but the product of hoping rates from $C_{i}$ to $C_{i+1}$ and from $C_{i+1}$ to $C_{i}$. 

We define vectors $\Psi_{i} \in \dC^{|V|}$ for $i \in [0, n]$ by $\Psi_{i}= (1/\sqrt{|C_{i}|})  \bm{1}_{i}$, where $\bm{1}_{i}$ is a vector in $\in \dC^{|V|}$ with $\bm{1}_{i}(v)=1$ for $v \in C_{i}$ and $\bm{1}_{i}(v)=0$ for $v \notin C_{i}$. In addition, we define the following subspace on $\dC^{|V|}$.
\[ \mathcal{A}=\mathrm{Span} \{ \Psi_{0}, \Psi_{1}, \cdots, \Psi_{n} \} .\]
Then, it holds that 
\[ \mathcal{A}^{\bot}= \left\{ f \in \dC^{|V|} \bigg{|} \sum_{v \in C_{i}} f(v) = 0, \, \forall i \in [0, n] \right\}. \]
Let $\mathcal{U}$ be the diagonal matrix of vertex degrees. We consider the following symmetric matrix:
\[
\tilde{T}=\mathcal{U}^{1/2}T\mathcal{U}^{-1/2}.
\]
We often use $\tilde{T}$ instead of $T$ since the symmetry of the matrix helps us to analyze its eigenvalues. We also denote the representation matrix of $\tilde{T}$ on a subspace $X \in \dC^{|V|}$ by $\tilde{T}|_{X}$. Then, the following holds.
\begin{lem}
For a generalized Bethe tree $B(d(0), d(1), \cdots, d(n-1))$, $\tilde{T}|_{\mathcal{A}}$ can be written by the following $(n+1) \times (n+1)$ tri-diagonal matrix in the basis of $\{ \Psi_{0}, \Psi_{1}, \cdots, \Psi_{n} \}$:
\begin{equation*}
\tilde{T}|_{\mathcal{A}}=\left(
	\begin{array}{ccccccc}
	0 & \sqrt{D_{0}} & & & & &   \\
	\sqrt{D_{0}} & 0 & \sqrt{D_{1}} & & & &  \\
	 & \sqrt{D_{1}} & 0 & \sqrt{D_{2}} & & &   \\	 
	 & & & \ddots & & &   \\
	 & & & & \ddots & &   \\
	 & & & & & \ddots &  \\
	 & & & & \sqrt{D_{n-2}} & 0 & \sqrt{D_{n-1}} \\
	 & & & & & \sqrt{D_{n-1}} & 0 \\
	 \end{array}
	\right).
\end{equation*} \label{rem of T}
\end{lem}
\begin{proof}
From the definition of the general Bethe tree, we have $|C_{0}|=1$, and $|C_{i}|=|C_{i-1}|d(i-1)$ for $i \in [1, n]$. Thus, it holds 
\begin{align*}
\tilde{T} \Psi_{0}&=\frac{1}{\sqrt{d(0)}} \frac{1}{\sqrt{d(1)+1}} \bm{1}_{1}=\sqrt{D_{0}} \Psi_{1}, 
\end{align*}
and for $i \in [1, n-1]$, 
\begin{align*}
\tilde{T} \Psi_{i}&=\frac{d(i-1)}{\sqrt{d(i-1)+1}} \frac{1}{\sqrt{d(i)+1}} \frac{1}{\sqrt{|C_{i}|}} \bm{1}_{i-1}+\frac{1}{\sqrt{d(i)+1}} \frac{1}{\sqrt{d(i+1)+1}} \frac{1}{\sqrt{|C_{i}|}} \bm{1}_{i+1}\\
&=\sqrt{D_{i}} \Psi_{i+1}+\sqrt{D_{i-1}} \Psi_{i-1},
\end{align*}
and
\begin{align*}
\tilde{T} \Psi_{n}&=\frac{d(n-1)}{\sqrt{d(n-1)+1}} \frac{1}{\sqrt{d(n-1)}} \bm{1}_{n-1}=\sqrt{D_{n-1}} \Psi_{n-1}.
\end{align*}
\end{proof}

\subsection{Fundamental analysis tools for Grover walk on generalized Bethe trees}
 
 To prove Theorem \ref{thm of bethe}, we provide some key tools for the Grover walk on the generalized Bethe trees. For a generalized Bethe tree $G=B(d(0), d(1), \cdots, d(n-1))$, let $P$ be a map from $V\backslash C_{0}$ to $V$, which maps $v \in C_{i}$ to its parent vertex, which is the unique vertex in $N(v)\cap C_{i-1}$, for $i \in [1, n]$. Also, we denote 
 \[ P^{i}(v)=\underbrace{P(P(\cdots P(P(v)))}_{i}, \]
 and define $P^{0}(v)=v$. For $v \in V$, we give the subset $N_{n}(v) \subset V(G)$ such that
\[ N_{n}(v)= \{ w \in C_{n}\,| \,{}^\exists j, P^{j}(w)=v \}, \]
which is the set of the descendants of $v$ in $C_{n}$. 

We define a sequence of polynomials of $\lambda \in \dC$ with the following inductive procedure:
\[
\begin{cases}
g_{0}(\lambda)= 1,  \\
g_{1}(\lambda)= \lambda,   \\ 
g_{i}(\lambda)=  \left( d(n-i+1)+1 \right) \lambda g_{i-1}(\lambda)-d(n-i+1) g_{i-2}(\lambda), &  i \in [2, n], \\ 
g_{n+1}(\lambda)= d(0)(\lambda g_{n}(\lambda)-g_{n-1}(\lambda)).  \\
\end{cases} 
\] 

Using these polynomials,  we give the following Lemma and Theorem to analyze the eigenvalues of the transition matrix of $B(d(1), d(2), \cdots, d(n-1))$.
We also define
\begin{equation}
\begin{cases}
p_{0}(\lambda)=1, \\
p_{i}(\lambda)=\mathrm{det}(\lambda I_{i}-\tilde{T}^{(i)}), & i \in [1, n+1],
\end{cases}
\label{def of p}
\end{equation}
where $I_{i}$ is $i \times i$ identity matrix and $\tilde{T}^{(i)}$ is the subprincipal matrix of $\tilde{T}|_{\mathcal{A}}$ with removing from the first row to the $(n+1-i)$-th row and from the first column to the $(n+1-i)$-th column for $i \in [1, n+1]$. Let $\tilde{g}_{i}(\lambda)$ be the monic polynomial of $g_{i}(\lambda)$ for $i \in [0, n+1]$. Then, the following holds. 
\begin{lem}
It holds \[ p_{i}(\lambda)=\tilde{g}_{i}(\lambda) \] for any $\lambda \in \dC$ and $i \in [0, n+1]$. \label{lem of p g}
\end{lem} 
\begin{proof}
Clearly, $\tilde{g}_0(\lambda) = p_0(\lambda) = 1$ and
$\tilde{g}_1(\lambda) = p_1(\lambda) = \lambda$,
so in order to complete the proof,
it is enough to show that $\{\tilde{g}_i\}$ and $\{p_i\}$ satisfy the same recurrence relation.
First, by cofactor expansion, we have
\[ p_i(\lambda) = \lambda p_{i-1}(\lambda) - D_{n-i+1} p_{i-2}(\lambda) \]
for $i \in [0,n+1]$. By the construction of $\{ g_{i} \}$, the coefficient of the maximum degree of $g_{i}(\lambda)$ is $\prod^{i-1}_{j=1} (d(n-j)+1)$, 
for $i \in [2, n]$ and that of $g_{n+1}(\lambda)$ is $d(0) \prod^{n-1}_{j=1} (d(n-j)+1)$ 
Therefore, 
\[ \tilde{g}_{i}(\lambda)=\frac{g_{i}(\lambda)}{ \prod^{i-1}_{j=1} (d(n-j)+1) } \]
for $i \in [2, n]$. We can check $\tilde{g}_i(\lambda)
= \lambda \tilde{g}_{i-1}(\lambda) - D_{n-i+1} \tilde{g}_{i-2}(\lambda)$
for $i \in \{2,3\}$ since $d(n) = 0$.
For $i \in [4,n]$,
\begin{align*}
\tilde{g}_i(\lambda)
&= \frac{(d(n-i+1)+1)\lambda g_{i-1}(\lambda)}{\prod_{j=1}^{i-1}(d(n-j)+1)} + 
     \frac{d(n-i+1)g_{i-2}(\lambda)}{\prod_{j=1}^{i-1}(d(n-j)+1)} \\
&= \frac{\lambda g_{i-1}(\lambda)}{\prod_{j=1}^{i-2}(d(n-j)+1)} +
     \frac{d(n-i+1)g_{i-2}(\lambda)}{(d(n-i+1)+1)(d(n-i+2)+1) \prod_{j=1}^{i-3}(d(n-j)+1)} \\
&= \lambda \tilde{g}_{i-1}(\lambda) - D_{n-i+1} \tilde{g}_{i-2}(\lambda).
\end{align*}
We also have $\tilde{g}_{n+1}(\lambda) = \lambda \tilde{g}_{n}(\lambda) - D_0 \tilde{g}_{n-1} (\lambda)$
since $g_{n+1}(\lambda) = d(0)(\lambda g_n(\lambda) - g_{n-1}(\lambda))$
and $D_{0} = 1/(d(1)+1)$,
so we can see that $p_i(\lambda) = \tilde{g}_i(\lambda)$ for any $\lambda \in \dC$ and $i \in [0,n+1]$.
\end{proof} Then, it immediately follows that the roots of $p_{n+1}(\lambda)$ are the eigenvalues of $\tilde{T}$ on $\mathcal{A}$. Indeed, the polynomials $\{p_{i} \}^{n}_{i=0}$ are the monic orthogonal polynomials with the finite Jacobi coefficients $\{ \sqrt{D_{0}}, \sqrt{D_{1}}, \cdots, \sqrt{D_{n-1}} \}$ and then, the eigenfunction for $\lambda$ with $p_{n+1}(\lambda)=0$ is written by 
\begin{equation}
{}^{T}(\mathcal{D}_{n}p_{n}(\lambda), \mathcal{D}_{n-1}p_{n-1}(\lambda), \cdots, \mathcal{D}_{1}p_{1}(\lambda), \mathcal{D}_{0}(\lambda) ),
\label{ef of a}
\end{equation}
where $\mathcal{D}_{i}=1/ \left( \prod^{i}_{j=0}\sqrt{D_{n-1-i}} \right)$ \cite{Obata}. Also, the set of levels with at least $2$ children is denoted by $\Omega$, that is, 
\[ \Omega:=\{ i \in [1, n]\,|\,d(n-i) \ge 2 \}. \]
Then, we can find the eigenvalues and eigenfunctions of $T$ on $\mathcal{A}^{\bot}$.  

\begin{lem}
For every $i \in \Omega$, the roots of $g_{i}(\lambda)$ are eigenvalues of $T$ on $\mathcal{A}^{\bot}$. Moreover, for a fixed vertex $v^{*} \in C_{n-i}$ and $N(v^{*}) \cap C_{n-i+1} = \{ v_{1}, v_{2}, \cdots v_{l} \} $, one of its eigenfunction $f \in \dC^{|V|}$ of the eigenvalue $\lambda$ is given as follows:\\  
For  $w \in N_{n}(v_{1})$, $j \in [0, i-1]$, 
\[
  f(P^{j}(w))=g_{j}(\lambda);
\]
for $w \in N_{n}(v_{2})$, $j \in [0, i-1]$, 
\[
    f(P^{j}(w))=-g_{j}(\lambda);
\]
for $w \in C_{n} \backslash \left( N_{n}(v_{1}) \cup N_{n}(v_{2}) \right)$, $j \in [0, i-1]$,
\[
    f(P^{j}(w))=0;
\]
for $w \in \cup^{n-i}_{j=0} C_{j}$, 
\[ f(w)=0. \]

\label{lem of g}
\end{lem}
\begin{proof}
We show that $(T f)(x)=\lambda f(x)$ for any $x \in V$. From the definition of $T$, for any $w \in N_{n}(v_{1})$, and $j \in [0, i-1]$, it holds that
\begin{align*}
 (T f)(P^{j}(w)) &= \sum_{u \sim P^{j}(w)} \frac{1}{\mathrm{deg}(P^{j}(w))} f(u)\\
&=\frac{1}{d(n-j)+1} \left( f(P^{j+1}(w)) + \sum_{\substack{u \sim P^{j}(w) \\ u \in C_{n-j+1}}} f(u) \right)\\
& = \frac{1}{d(n-j)+1} \left( g_{j+1}(\lambda) + d(n-j) g_{j-1}(\lambda) \right)\\
& = \lambda g_{j}(\lambda)\\
& = \lambda f(P^{j}(w)).
\end{align*}
Similarly,  for any $w \in N_{n}(v_{2})$, and every $j \in [0, i-1]$, it holds that
\[ (T f)(P^{j}(w)) = -\lambda g_{j}(\lambda) = \lambda f(P^{j}(w)). \]
Furthermore, it is trivial that for any $w \in C_{n} \backslash \left( N_{n}(v_{1}) \cup N_{n}(v_{2}) \right)$, and every $j \in [0, i-1]$, 
\[ (T f)(P^{j}(w)) =0=\lambda f(P^{j}(w)). \]
Thus, we can say that $(T f)(x) = \lambda f(x)$ for $x \in \cup^{n}_{j=n-i+1} C_{j}$. In order to complete the proof, it is sufficient to show that $(T f)(v^{*})=\lambda f(v^{*})$ since it trivially holds for any $u \in \cup^{n-i}_{j=0}C_{j}\backslash \{v^{*} \}$. Indeed, we can say that
\begin{align*}
(T f)(v^{*}) = \frac{1}{\mathrm{deg}(v^{*})} \left( g_{i-1}(\lambda) - g_{i-1}(\lambda) \right)= 0 = \lambda f(v^{*}),
\end{align*}
from the definition of $f$. Hence, it holds that $T f =\lambda f$. Also, we have $\sum_{u \in C_{i}}f(u)=0$ for any $i \in [0, n]$ from the definition of $f$. Then $f \in \mathcal{A}^{\bot}$. 
\end{proof}

Therefore, we can obtain $\sigma(\tilde{T}|_{\mathcal{A}^{\bot}}) \supset \bigcup_{i \in \Omega} \{ \lambda \in \dC\,|\, g_{i}(\lambda)=0 \}$. Moreover, we can show its equality. 
\begin{lem}
\[ 
\sigma(\tilde{T}|_{\mathcal{A}^{\bot}}) = \bigcup_{i \in \Omega} \{ \lambda \in \dC\,|\, g_{i}(\lambda)=0 \}.
\] \label{prop of abot}
\end{lem}  
\begin{proof}
We show that the total of the linearly independent eigenfunctions on the above argument coincides with $\dim \mathcal{A^{\bot}}=|V|-(n+1)$. From the definition of the generalized Bethe tree, we have $|C_{0}|=1$ and 
\begin{equation}
|C_{i}|=\prod^{i-1}_{j=0}d(j) \label{num of ci}
\end{equation}
for $i \in [1, n-1]$. Thus,  
\[ |V| = \sum^{n}_{i=0} |C_{i}| = |C_{0}|+\sum^{n-1}_{i=0}|C_{n-i}|=1+\sum^{n-1}_{i=0} \prod^{n-i-1}_{j=0}d(j). \]
So it holds that $\dim \mathcal{A}^{\bot}=  \sum^{n-1}_{i=0} \prod^{n-i-1}_{j=0}d(j)-n$.  We fix an integer $i \in \Omega$ and let $\lambda$ be a root of $g_{i}$. By Theorem \ref{lem of g}, there are $d(n-i)-1$, which is the number of combinations of the children $v_{1}$ and the other vertices,  linearly independent eigenfunctions corresponding to $\lambda$ for every $v \in C_{n-i}$. So we can make $|C_{n-i}|\left( d(n-i)-1 \right)$ eigenfunctions corresponding to $\lambda$. The roots of $g_{i}$ are the eigenvalues of a tri-diagonal matrix and all the eigenvalues of the tri-diagonal matrices are distinct \cite{tri}.  Hence, the total of the linearly independent eigenfunctions corresponding to the roots of $g_{i}$ is $i |C_{n-i}|\left( d(n-i)-1 \right)$. Here, we put $F_{i}:=i |C_{n-i}|\left( d(n-i)-1 \right)$ for $i \in [1, n]$. By (\ref{num of ci}), it follows
\[
F_{i} =  i \left( \prod^{n-i}_{j=0}d(j)-\prod^{n-i-1}_{j=0}d(j) \right),
\]
for $i \in [1, n-1]$ and $F_{n}=n(d(0)-1)$. Then, the total of the linearly independent eigenfunctions for the roots of every $g_{i}$ for $i \in \Omega$ is $\sum^{n}_{i=1} F_{i}$ because $F_{i}=0$ for $i$ with $d(n-i)=1$. Thus, we can get 
\begin{align*}
\sum_{i=1}^{n} F_i 
&= \sum_{i=1}^{n-1} i \left( \prod_{j = 0}^{n-i} d(j) - \prod_{j = 0}^{n-i-1} d(j) \right) + F_n \\
&= \sum_{i=1}^{n-1} \left( i \prod_{j = 0}^{n-i} d(j) - i \prod_{j = 0}^{n-i-1} d(j) \right) + F_n \\
&= \sum_{i=1}^{n-1} \prod_{j=0}^{n-i} d(j) - (n-1)d(0) + n(d(0) - 1) \\
&= \sum_{i=1}^{n-1} \prod_{j=0}^{n-i} d(j) + d(0) - n \\
&= \sum_{i=1}^{n} \prod_{j=0}^{n-i} d(j) - n,
\end{align*}
which equals $\dim \mathcal{A}^{\bot}$. 
\end{proof}
Let $f^{(\lambda)}_{v^{*} : i} \in \dC^{|V|}$ be the eigenfunction mentioned in Lemma \ref{lem of g} for $i \in \Omega$, $v^{*} \in C_{n-i}$, and $\lambda$ with $p_{i}(\lambda)=0$. Also, for $\lambda$ with $p_{n+1}(\lambda)=0$, we denote the vector in (\ref{ef of a}) by $f^{(\lambda)}_{n+1}$. 
Then, by Lemmas \ref{lem of g}, \ref{prop of abot} and the equation (\ref{ef of a}), we can summarize the above statements as follows: 
\begin{thm}
\[ \sigma(\tilde{T}|_{\mathcal{A}^{\bot}})=\bigcup_{i \in \Omega} \{ \lambda \in \dC \,|\,p_{i}(\lambda)=0 \},\]
\[
\sigma(\tilde{T}|_{\mathcal{A}})=\{ \lambda \in \dC \,|\,p_{n+1}(\lambda)=0 \},
\]
\[ \mathrm{ker}(\lambda I-\tilde{T}_{\mathcal{A}^{\bot}})=\mathrm{Span}\{ f^{(\lambda)}_{v^{*} : i} \,|\, i \in \Omega, v^{*} \in C_{n-i}, p_{i}(\lambda)=0 \}, \]
\[ \mathrm{ker}(\lambda I - \tilde{T}|_{\mathcal{A}})=\mathrm{Span}\{f^{(\lambda)}_{n+1}\,|\, p_{n+1}(\lambda)=0 \}. \]
\end{thm}

Before proving our main Theorem, recall the two kinds of the Chebyshev polynomials and give a useful Lemma. Let $T_i$ and $U_i$ be the Chebyshev polynomial of the first kind of degree $i$
and the one of the second kind of degree $i$, respectively.
Then, it is well-known that $T_{0}(\lambda)=1, T_{1}(\lambda)=\lambda$, $U_{-1}(\lambda)=0, U_{0}(\lambda)=1, U_{1}(\lambda)=2\lambda$, and for $i \ge 2$, these polynomials are defined by the same recurrence relation such that
\[ L_{i}(\lambda)=2\lambda L_{i-1}(\lambda)-L_{i-2}(\lambda), \]
where $L_{i}=T_{i}$, or $U_{i}$. By the above relation, it follows that 
\begin{align}
T_i(\cos \theta) &= \cos (i \theta) \label{T1} \\
\intertext{and}
U_i(\cos \theta) &= \frac{\sin ((i+1)\theta)}{\sin \theta}. \label{U1}
\end{align}
Then, we will provide the following Lemma. 
\begin{lem}
For $z \in \dC$ with $|z| = 1$ and a positive integer $i$,
the followings hold.
	\begin{description}
	\item[(i)] $(2z)^i T_i \left( \dfrac{z+z^{-1}}{2} \right) = 2^{i-1}(z^{2i} + 1)$.
	\item[(ii)] $\displaystyle (2z)^i U_i \left( \dfrac{z+z^{-1}}{2} \right)= 2^i \sum_{j=0}^i z^{2i-2j}$.
	\end{description}
	\label{lem of cheby}
\end{lem}
It can be proved by putting $z=e^{\sqrt{-1}\theta}$, and it will be useful to show our main result. 
\section{Proof of main result}
In this section, we prove our main result with some tools. 
\begin{proof}[Proof of Theorem \ref{thm of bethe}] For the generalized Bethe tree $B(d(0), d(1), \cdots, d(n))$, let $k_{1}, k_{2}, \cdots, k_{l}$ be positive integers such that for $K_{i}=\sum^{i}_{j=1} k_{j}$,  $C_{n-K_{i}}$ is the only level with $d(n-K_{i}) \ge 2$ for $i \in [1, l]$ as is seen in Figure \ref{fig bethe5}. Thus, $l$ is the number of levels having at least $2$ children. This value will be important in this proof. We define $k_{l+1}=n-K_{l}$, then $K_{l+1}=k_{l+1}+K_{l}=n$. Then, we can regard each level as {\it age}, where $n$-th level is $0$ age. So for $i \in [1, l]$, $K_{i}$ is an age having at least $2$ children. Also, we put $d_{i}=d(n-K_{i})$, which is the number of children of $K_{i}$ age. By Lemmas \ref{lem of p g} and \ref{lem of g}, the roots of $p_{K_{i}}(\lambda)$ are the eigenvalue of $T$ on $\mathcal{A}^{\bot}$ for $i \in [1, l]$, and those of $p_{K_{(l+1)}+1}$ are the eigenvalues of $T$ on $\mathcal{A}$. So we find the polynomials $p_{K_{i}}(\lambda)$, and check whether $(2z)^{K_{i}}p_{K_{i}}\left( (z+z^{-1})/2 \right)$ is represented by a product of some cyclotomic polynomials for every $i \in [1, l+1]$. For simplicity, we write a polynomial $f(\lambda)$ as $f$ in this section.
\begin{figure}[h]
		\centering
		\scalebox{1.0}{
		\includegraphics[scale=0.9]{bethe-figfig.3}
		}
		\caption{Setting}
		\label{fig bethe5}
		\end{figure}\\
Here, we introduce the following Claims.
\begin{claim}
Let $p_{i}$ be defined in (\ref{def of p}). For $i \in [1, K_{1}]$, we have
\begin{equation}
p_{i}=\frac{1}{2^{i-1}}T_{i}.  
\end{equation}
\label{eq: till pk1}
\end{claim}
\begin{proof}
If $K_{1}=1, 2$, then it clearly holds because $p_{1}=\lambda=T_{1}$, and $p_{2}=\lambda p_{1}-D_{n-1}p_{0}=\lambda^{2}-\frac{1}{2}=\frac{1}{2} T_{2}$. So we assume $K_{1} \ge 3$. We prove the above equality by the induction on $i$. Indeed, $D_{n-i+1}=1/4$ for $i \in [3, K_{1}]$, so it holds 
\begin{align*}
p_{i}&=\lambda p_{i-1}-D_{n-i+1} p_{i-2}\\
&=\lambda \left\{ \frac{1}{2^{i-2}} T_{i-1} \right\}-\frac{1}{4}  \left\{ \frac{1}{2^{i-3}} T_{i-2} \right\}\\
&=\frac{1}{2^{i-1}} \left( 2\lambda T_{i-1}-T_{i-2} \right)\\
&=\frac{1}{2^{i-1}} T_{i}.
\end{align*} 
\end{proof}
\begin{claim}
For $j \in [2, l+1]$ with $k_{j} \ge 2$, we have
\begin{equation}
p_{K_{(j-1)}+i}=\frac{1}{2^{i-2}} U_{i-2} p_{K_{(j-1)}+2}-\frac{1}{2^{i-1}} U_{i-3} p_{K_{(j-1)}+1}
\end{equation}
for $i \in [2, k_{j}]$. 
\label{eq: till pk2}
\end{claim}
\begin{proof}
Indeed, we can also prove the above equality by the induction on $i$. It is easy to confirm that this Claim holds for $i = 2, 3$. So we assume that $i \ge 4$. Then,
\begin{align*}
p_{K_{(j-1)}+i}&=\lambda p_{K_{(j-1)}+i-1}-D_{n-K_{(j-1)}-i+1} p_{K_{(j-1)}+i-2}\\
&=\lambda p_{K_{(j-1)}+i-1}-\frac{1}{4} p_{K_{(j-1)}+i-2}\\
&=\lambda \left\{  \frac{1}{2^{i-3}} U_{i-3} p_{K_{(j-1)}+2}-\frac{1}{2^{i-2}} U_{i-4} p_{K_{(j-1)}+1} \right\}-\frac{1}{4} \left\{ \frac{1}{2^{i-4}} U_{i-4} p_{K_{(j-1)}+2}-\frac{1}{2^{i-3}} U_{i-5} p_{K_{(j-1)}+1} \right\}\\
&= \frac{1}{2^{i-2}} ( 2 \lambda U_{i-3} - U_{i-4} )p_{K_{(j-1)}+2}
- \frac{1}{2^{i-1}}( 2 \lambda U_{i-4} - U_{i-5}) p_{K_{(j-1)}+1} \\
&= \frac{1}{2^{i-2}} U_{i-2} p_{K_{(j-1)}+2} - \frac{1}{2^{i-1}} U_{i-3} p_{K_{(j-1)}+1}.
\end{align*}
\end{proof}

First, we suppose $l=0$, then, the generalized Bethe tree is nothing but the path graph $P_{n+1}$. So all the eigenvalues of $T$ coincide with the roots of $p_{K_{1}+1}=\frac{1}{2^{n-1}}\left\{ \lambda T_{n}-T_{n-1} \right\}$ by Claim \ref{eq: till pk1}. Then, we have
\[ (2z)^{K_{1}+1} p_{K_{1}+1}\left( \frac{z+z^{-1}}{2} \right)= z^{2n+2}-z^{2n}-z^{2}+1=(z^{2}-1)(z^{2n}-1)\]
 by (i) on Lemma \ref{lem of cheby}, and its roots satisfy $z^{2n}=1$. Hence, the path graph $P_{n+1}$ induces a $2n$-periodic Grover walk. 

Next, we suppose $l=1$. Then, the generalized Bethe tree is  
\[ B(\underbrace{1, 1, \cdots, 1}_{k_{2}}, d_{1}, \underbrace{1, 1, \cdots, 1}_{k_{1}-1}). \] 
If $k_{2}=0$ and $k_{1}=1$, or $k_{2}=1$ and $k_{1}=1$, then these graphs are nothing but star graphs, which induce a $4$-periodic Grover walk \cite{Yoshie}. Hence, we do not have to consider these cases. 
By claim \ref{eq: till pk1}, we can obtain 
\[ p_{K_{1}}=\frac{1}{2^{K_{1}-1}}T_{K_{1}}. \]
By (i) on Lemma \ref{lem of cheby}, we have
\begin{equation}
(2z)^{K_{1}} p_{K_{1}}\left( \frac{z+z^{-1}}{2} \right)= z^{2K_{1}}+1. \label{gk1}
\end{equation}
Thus, the eigenvalues of $\tilde{T}$ on $\mathcal{A^{\bot}}$ satisfy the condition on Lemma \ref{lem of cyclo} and the roots satisfy 
\begin{equation}
z^{4k_{1}}=1. \label{sol of gk1}
\end{equation}
We have to find the polynomial $p_{n+1}=p_{K_{2}+1}$ to analyze the eigenvalues of $\tilde{T}$ on $\mathcal{A}$ in this case. Then, we have
\begin{equation*}
p_{K_{1}+1}=
	\begin{cases}
	\lambda p_{1}-\frac{d_{1}}{d_{1}+1} p_{0} &\text{if $k_{1}=1$,}\\
	\lambda p_{K_{1}}-\frac{d_{1}}{2(d_{1}+1)} p_{K_{1}-1} & \text{if $k_{1} \ge 2$,}
	\end{cases}
\end{equation*}
and
\begin{equation*}
p_{K_{1}+2}=
	\begin{cases}
	\lambda p_{K_{1}+1}-\frac{1}{d_{1}+1} p_{K_{1}} &\text{if $k_{2}=1$,}\\
	\lambda p_{K_{1}+1}-\frac{1}{2(d_{1}+1)} p_{K_{1}} & \text{if $k_{2} \ge 2$.}
	\end{cases}
\end{equation*}
Indeed, for any $k_{1} \in \dN$, we can calculate that  
\begin{equation}
(2z)^{K_{1}+1} p_{K_{1}+1}\left( \frac{z+z^{-1}}{2} \right)=z^{2K_{1}+2}+\frac{1-d_{1}}{d_{1}+1}z^{2K_{1}}+\frac{1-d_{1}}{d_{1}+1}z^{2}+1. \label{pk1+1}
\end{equation}
Furthermore, we can also calculate that if $k_{2}=1$, then
\begin{equation}
(2z)^{K_{1}+2} p_{K_{1}+2}\left( \frac{z+z^{-1}}{2} \right)=
z^{2K_{1}+4}-\frac{2}{d_{1}+1}z^{2K_{1}+2}+\frac{1-d_{1}}{d_{1}+1}z^{2K_{1}}+\frac{1-d_{1}}{d_{1}+1}z^{4}-\frac{2}{d_{1}+1}z^{2}+1 \label{pk1+2.1}
\end{equation}
and if $k_{2} \ge 2$, then
\begin{equation}
(2z)^{K_{1}+2} p_{K_{1}+2}\left( \frac{z+z^{-1}}{2} \right)=z^{2K_{1}+4}+\frac{1-d_{1}}{d_{1}+1}z^{2K_{1}}+\frac{1-d_{1}}{d_{1}+1}z^{4}+1 \label{pk1+2}
\end{equation}
by (i) on Lemma \ref{lem of cheby}. 
Then, it follows that
\begin{equation*}
p_{K_{2}+1}=
	\begin{cases}
	p_{K_{1}+2} & \text{if $k_{2}=1$,}\\
	\lambda p_{K_{1}+2}-\frac{1}{2} p_{K_{1}+1} & \text{if $k_{2}=2$,}\\
         \frac{1}{2^{k_{2}-2}}\left( \lambda U_{k_{2}-2}-U_{k_{2}-3} \right)p_{K_{1}+2}-\frac{1}{2^{k_{2}-1}}\left( \lambda U_{k_{2}-3}-U_{k_{2}-4} \right) p_{K_{1}+1}. & \text{if $k_{2} \ge 3$,}\\
	\end{cases}
\end{equation*} 
by Claim \ref{eq: till pk2}. 
Using (\ref{pk1+1}), (\ref{pk1+2.1}), (\ref{pk1+2}) and (ii) on Lemma \ref{lem of cheby}, we can get 
\begin{align*}
(2z)^{K_{2}+1}p_{K_{2}+1}\left( \frac{z+z^{-1}}{2} \right)&=z^{2K_{2}+2}-z^{2K_{2}}+\frac{1-d_{1}}{d_{1}+1}\left( z^{2k_{2}+2}-z^{2k_{2}}-z^{2K_{1}+2}+z^{2K_{1}} \right) -z^{2}+1.
\end{align*}
Unless $k_{1}=K_{1}=k_{2}$, the above polynomial is not an integer polynomial unless, which implies that it cannot be represented by a product of some cyclotomic polynomials. Thus, the Bethe tree is nothing $k_{1}$-subdivision of $ST_{d_{1}+2}$. Then, the polynomial can be rewritten as
\[ z^{4k_{1}+2}-z^{4k_{1}}-z^{2}+1=(z^{2}-1)(z^{4k_{1}}-1) \]
and the roots satisfy $z^{4k_{1}}=1$. So the period is $4k_{1}$.

Next, we suppose $l \ge 2$. Then, the following holds
\begin{claim}
It should hold that $k_{1}=k_{2}$ and $d_{1}=3$ for $l \ge 2$. \label{claim of d}
\end{claim}
\begin{proof}
By Claim \ref{eq: till pk2}, we can obtain
\begin{equation*}
p_{K_{2}}=
	\begin{cases}
	p_{K_{1}+1} & \text{if $k_{2}=1$,}\\
	\frac{1}{2^{k_{2}-2}}U_{k_{2}-2}p_{K_{1}+2}-\frac{1}{2^{k_{2}-1}}U_{k_{2}-3}p_{K_{1}+1} & \text{if $k_{2} \ge 2$.}
	\end{cases}
\end{equation*}
Indeed, for any $k_{1}, k_{2} \in \dN$, it follows
\begin{equation}
(2z)^{K_{2}} p_{K_{2}}\left(\frac{z+z^{-1}}{2} \right)=z^{2K_{2}}+\frac{1-d_{1}}{d_{1}+1}z^{2k_{1}}+\frac{1-d_{1}}{d_{1}+1}z^{2k_{2}}+1 \label{eq of gk2}
\end{equation}
by (\ref{pk1+1}), (\ref{pk1+2.1}), (\ref{pk1+2}), and (ii) of Lemma \ref{lem of cheby}. Thus, (\ref{eq of gk2}) is not an integer polynomial unless $k_{1}=k_{2}$ and $d_{1}=3$. 
\end{proof}
Then, (\ref{eq of gk2}) turns to 
\begin{equation}
z^{4k_{1}}-z^{2k_{1}}+1, \label{gk2}
\end{equation}
and its roots satisfy 
\begin{equation}
z^{12k_{1}}=1. \label{sol of gk2}
\end{equation}
Furthermore, the following also holds to induce a periodic Grover walk in this case.
\begin{claim}
For $l \ge 2$ case, $l$ is exactly $2$. \label{claim of l}
\end{claim}
\begin{proof}
We suppose $l \ge 3$ and lead a contradiction. If $l \ge 3$, then the roots of $p_{K_{3}}$ should be the eigenvalue of $T$ on $\mathcal{A^{\bot}}$. Then, we also have
\begin{equation*}
p_{K_{3}}=
	\begin{cases}
	p_{K_{2}+1} & \text{if $k_{3}=1$,}\\
	\frac{1}{2^{k_{3}-2}}U_{k_{3}-2}p_{K_{2}+2}-\frac{1}{2^{k_{3}-1}}U_{k_{3}-3}p_{K_{2}+1} & \text{if $k_{3} \ge 2$.}
	\end{cases}
\end{equation*}
by Claim \ref{eq: till pk2}. Also, we can set $k_{1}=k_{2}$, and $d_{1}=3$ from Claim \ref{eq: till pk2}. Then, it follows that
\begin{align*}
p_{K_{2}+2}&=\lambda p_{K_{2}+1}-\frac{1}{2(d_{2}+1)} p_{K_{2}} \\
p_{K_{2}+1}&=
	\begin{cases}
	\lambda p_{K_{2}}-\frac{d_{2}}{4(d_{2}+1)} p_{K_{2}-1} &\text{if $k_{1}=k_{2}=1$, }\\
	\lambda p_{K_{2}}-\frac{d_{2}}{2(d_{2}+1)} p_{K_{2}-1} &\text{if $k_{1}=k_{2} \ge 2$.}\\
	\end{cases}\\
p_{K_{2}-1}&=
	\begin{cases}
	p_{K_{1}} &\text{if $k_{1}=k_{2}=1$,}\\
	p_{K_{1}+1} &\text{if $k_{1}=k_{2}=2$,}\\
	\frac{1}{2^{k_{2}-3}}U_{k_{2}-3}p_{K_{1}+2}-\frac{1}{2^{k_{2}-2}}U_{k_{2}-4}p_{K_{1}+1} &\text{if $k_{1}=k_{2} \ge 3$.}
	\end{cases}
\end{align*}
Combining (\ref{gk1}), (\ref{pk1+1}), (\ref{pk1+2.1}), (\ref{pk1+2}), and (\ref{gk2}), we have
\begin{equation}
(2z)^{K_{2}-1} p_{K_{2}-1}\left(\frac{z+z^{-1}}{2} \right)=z^{4k_{1}-2}+\frac{1}{2}z^{2k_{1}}-\frac{1}{2}z^{2k_{1}-2}+1,
\end{equation}
\begin{equation} 
(2z)^{K_{2}+1} p_{K_{2}+1}\left(\frac{z+z^{-1}}{2} \right)=z^{4k_{1}+2}+\frac{1-d_{2}}{d_{2}+1}z^{4k_{1}}-\frac{1}{d_{2}+1}z^{2k_{1}+2}-\frac{1}{d_{2}+1}z^{2k_{1}}+\frac{1-d_{2}}{d_{2}+1}z^{2}+1, 
\end{equation}
\begin{equation}
(2z)^{K_{2}+2} p_{K_{2}+2}\left(\frac{z+z^{-1}}{2} \right)=z^{4k_{1}+4}+\frac{1-d_{2}}{d_{2}+1}z^{4k_{1}}-\frac{1}{d_{2}+1}z^{2k_{1}+4}-\frac{1}{d_{2}+1}z^{2k_{1}}+\frac{1-d_{2}}{d_{2}+1}z^{4}+1.
\end{equation}
Thus, for any $k_{1}, k_{3} \in \dN$, it follows
\begin{equation}
(2z)^{K_{3}} p_{K_{3}}\left(\frac{z+z^{-1}}{2} \right)=z^{2K_{3}}+\frac{1-d_{2}}{d_{2}+1}z^{4k_{1}}-\frac{1}{d_{2}+1}z^{2k_{1}}-\frac{1}{d_{2}+1}z^{2k_{1}+2k_{3}}+\frac{1-d_{2}}{d_{2}+1}z^{2k_{3}}+1. \label{eq of gk3}
\end{equation}
The above polynomial cannot be an integer polynomial for any $d_{2} \in \dN_{\ge 2}$, and $k_{1}, k_{3} \in \dN$. In other words, the roots of $p_{K_{3}}$ cannot be roots of the real part of the unity for any $d_{2} \in \dN_{\ge 2}$, and $k_{1}, k_{3} \in \dN$. Therefore, we have to set $l < 3$ so that the roots of $p_{K_{3}}(\lambda)$ are not eigenvalues of $T$. 
\end{proof}
Claims \ref{claim of d}, \ref{claim of l} are necessary conditions only induced by $\mathcal{A}^{\bot}$. Therefore, the candidates of generalized Bethe trees which can induce a periodic Grover walk are written by
\begin{equation}
B(\underbrace{1, 1, \cdots, 1}_{k_{3}}, d_{2}, \underbrace{1, 1, \cdots, 1}_{k_{1}-1}, 3, \underbrace{1, 1, \cdots, 1}_{k_{1}-1}), \label{bethe 1} 
\end{equation}
or,
\begin{equation}
B(d_{2}, \underbrace{1, 1, \cdots, 1}_{k_{1}-1}, 3, \underbrace{1, 1, \cdots, 1}_{k_{1}-1}). \label{bethe 2}
\end{equation}

For (\ref{bethe 1}), we consider the eigenvalues of $T$ on $\mathcal{A}$ to determine the parameters $k_{3} \in \dN$ and $d_{2} \in \dN_{\ge 2}$ of (\ref{bethe 1}). We have to find $p_{n+1}=p_{K_{3}+1}$ and analyze its roots since these are coincidence to the eigenvalues of $T$ on $\mathcal{A}$. Then, we have
\begin{equation*}
p_{K_{3}+1}=
	\begin{cases}
	\lambda p_{K_{3}}-\frac{1}{d_{2}+1} p_{K_{2}} & \text{if $k_{3}=1$,}\\
	\lambda p_{K_{3}}-\frac{1}{2} p_{K_{3}-1} & \text{if $k_{3} \ge 2$.}
	\end{cases}
\end{equation*}
For $k_{1}, k_{3} \in \dN$ the polynomials $(2z)^{K_{3}} p_{K_{3}}\left( (z+z^{-1})/2 \right), (2z)^{K_{2}} p_{K_{2}}\left( (z+z^{-1})/2 \right)$ are of the forms (\ref{eq of gk3}), (\ref{gk2}), respectively. Also, we can calculate 
\begin{equation*}
\begin{split}
 (2z)^{K_{3}-1} p_{K_{3}-1}\left( \frac{z+z^{-1}}{2} \right)=z^{2(K_{3}-1)}+\frac{1-d_{2}}{d_{2}+1}z^{4k_{1}}-\frac{1}{d_{2}+1}z^{2k_{1}}\\
 -\frac{1}{d_{2}+1}z^{2k_{1}+2(k_{3}-1)}+\frac{1-d_{2}}{d_{2}+1}z^{2(k_{3}-1)}+1. 
\end{split}
\end{equation*}
Indeed, for any $k_{3} \in \dN$, we can also calculate that $(2z)^{K_{3}+1}p_{K_{3}+1}\left((z+z^{-1})/2 \right)$ as follows:  
\begin{equation*}
\begin{split}
 (2z)^{K_{3}+1}p_{K_{3}+1}\left(\frac{z+z^{-1}}{2} \right)=& z^{2K_{3}+2}-\frac{1}{d_{2}+1} \left(z^{2k_{1}+2k_{3}+2}-z^{2k_{1}+2k_{3}}-z^{2k_{1}+2}+z^{2k_{1}} \right)\\
 &-\frac{1-d_{2}}{d_{2}+1}\left(z^{4k_{1}+2}-z^{4k_{1}}-z^{2k_{3}+2}+z^{2k_{3}} \right)-z^{2K_{3}}-z^{2}+1. 
 \end{split}
 \end{equation*}
The above polynomial also cannot be an integer polynomial unless $k_{1}=k_{3}$ and $d_{2}=2$. Then, the right hand side turns to  
\[ z^{6k_{1}+2}-z^{6k_{1}}-z^{2}+1=(z+1)(z-1)(z^{6k_{1}}-1), \]
whose roots satisfy
\begin{equation}
z^{6k_{1}}=1. \label{sol of a1}
\end{equation} 
Thus, it should hold $k_{3}=k_{1}$ and $d_{2}=2$. Put $k_{1}=k$. Then, the generalized Bethe tree is $B(\underbrace{1, 1, \cdots, 1}_{k}, 2, \underbrace{1, 1, \cdots, 1}_{k-1}, 3, \underbrace{1, 1, \cdots, 1}_{k-1})$, which is nothing but $S_{k}(B(1, 2, 3))$. Also, it follows that its period is $12k$ from (\ref{sol of gk1}), (\ref{sol of gk2}), (\ref{sol of a1}). 

For (\ref{bethe 2}), we find $p_{n+1}=p_{K_{2}+1}$ and analyze its roots. We also have
\begin{equation*}
p_{K_{2}+1}=
	\begin{cases}
	\lambda p_{K_{2}}-\frac{1}{4} p_{K_{1}} & \text{if $k_{1}=1$,}\\
	\lambda p_{K_{3}}-\frac{1}{2} p_{K_{2}-1} & \text{if $k_{1} \ge 2$.}
	\end{cases}
\end{equation*}
The polynomial $p_{K_{1}}$ is of the form (\ref{gk1}) and we also have 
\[ (2z)^{K_{2}-1} p_{K_{2}-1} \left( \frac{z+z^{-1}}{2} \right)=z^{4(k_{1}-1)}-z^{2(k_{1}-1)}+1 \]
for $k_{1} \in \dN$. Thus, for any $k_{1} \in \dN$, then we can calculate $(2z)^{K_{2}+1}p_{K_{2}+1}\left((z+z^{-1})/2 \right)$ as follows:  
\[ z^{4k_{1}+2}-z^{4k_{1}}-z^{2}+1=(z+1)(z-1)(z^{4k_{1}}-1), \]
whose roots satisfy 
\begin{equation}
z^{4k_{1}}=1. \label{sol of a2}
\end{equation} 
Put $k_{1}=k$ and $d_{2}=s$. Then, we can obtain that for any $k, s \in \dN$, $B(s, \underbrace{1, 1, \cdots, 1}_{k-1}, 3, \underbrace{1, 1, \cdots, 1}_{k-1})$, which is nothing but $S_{k}(B(s, 3))$, only induces a $12k$-periodic Grover walk by (\ref{sol of gk1}), (\ref{sol of gk2}), (\ref{sol of a2}) on the case (\ref{bethe 2}). Therefore, the proof is completed. 
\end{proof}

\section{Summary and discussions}
In this paper, we completely classified the generalized Bethe trees which induce a periodic Grover walk. 
These graphs are the path graphs and subdivision of the star graphs, $B(1, 2, 3), B(s, 3)$ for $s \in \dN_{\ge 2}$. Remark that the generalized Bethe trees are tree graphs with an equitable partition. So, it is natural to consider the periodicity of the Grover walk on other graphs with an equitable partition, e.g. spider net graphs, distance-regular graphs \cite{Kempe}, \cite{KOS}, \cite{Amb}. Since the spectrum of such graphs can be analyzed with orthogonal polynomials induced by a probability measure on not only real line but also the unit circle in the complex plane \cite{Cantero}, the method used in this paper can be useful to consider its periodicity. The Grover walk is given by a unitary operator called the Grover transfer matrix. If we use another unitary operator as an evolution operator, we can determine another quantum walk. So it is also our future work to consider the periodicity of quantum walks driven by different unitary operator from the Grover transfer matrix on graphs with an equitable partition. 

\section*{Acknowledgments}
E.S. acknowledges financial support from the Grant-in-Aid for Young Scientists (B) and of Scientific Research (B) Japan Society for the Promotion of Science (Grant No. 16K17637, No. 16K03939). T.T.  was supported by JSPS KAKENHI (Grant No.  16K05263).



\begin{thebibliography}{9}
  
  \bibitem{Amb}
	A.\,Ambainis,
	\newblock Quantum walk algorithm for element distinctness,
	\newblock SIAM Journal on Computing, {\bf 37}, pp. 210--239, 
	\newblock (2007)  
	
\bibitem{Cantero}
	M.\,J.\,Cantero, L.\,Moral L.\,Velazquez,
	\newblock Five-diagonal matrices and zeros of orthogonal polynomials on the unit circle,
	\newblock Linear Algebra and its Applications, {\bf 362}, pp. 29--56,
	\newblock (2003)	
	
\bibitem{Chisa11}
	K.\,Chisaki, N.\,Konno, E.\,Segawa, Y.\,Shikano,
	\newblock Crossovers induced by discrete-time quantum walks,
	\newblock Quantum Information and Computation, {\bf 11}, pp. 741--760, 
	\newblock (2011)	
	
\bibitem{Emms06}
  D.\,M.\,Emms, E.\,R.\,Hancock, S.\,Severini, R.\,C.\,Wilson,
  \newblock A matrix representation of graphs and its spectrum as a graph invariant, 
  \newblock Electronic Journal Combinatorics, {\bf 13},
  \newblock R34\, (2006)	
  
  \bibitem{RGodsil}
	C.\,Godsil, G.\,Royle,
	\newblock Algebraic Graph Theory, 
	\newblock Springer New York,
	\newblock (2001)
	
\bibitem{Godsil}
	C.\,Godsil,
	\newblock State transfer on graphs,
	\newblock Discrete Mathematics, {\bf 312}, pp. 129--147,
	\newblock (2011)
	
\bibitem{Grover}
	L.\,Grover,
	\newblock A fast quantum mechanical algorithm for database search, 
	\newblock Proceedings of the 28th Annual ACM Symposium on the Theory of Computing, pp. 212--219,
	\newblock (1996)
	
\bibitem{Gudder}
	S.\,Gudder,
	\newblock Quantum Probability,
	\newblock Academic Press Inc., CA, 
	\newblock (1998)	
	
\bibitem{HSegawa13}
	 Yu.\,Higuchi, N.\,Konno, I.\,Sato, E.\,Segawa,
	 \newblock A note on the discrete-time evolutions of quantum walk on a graph,
	 \newblock Journal of Math-for-Industry, {\bf 5}, pp. 103--109,
	 \newblock (2013)
	 
  \bibitem{Higuchi13}
  Yu.\,Higuchi, N.\,Konno, I.\,Sato, E.\,Segawa,
  \newblock Periodicity of the discrete-time quantum walk on finite graph,
  \newblock Interdisciplinary Information Sciences, {\bf 23}, pp. 75--86,
  \newblock (2017)
  
  \bibitem{HSegawa17}
	Yu.\,Higuchi, E.\,Segawa,
	\newblock  Quantum walks induced by Dirichlet random walks on infinite trees,
	\newblock accepted to publication to Journal of Physics A: Mathematical and Theoretical. 
	\newblock (2017)  
	
\bibitem{Obata}
	A.\,Hora, N.\,Obata,
	\newblock  Quantum Probability and Spectral Analysis of Graphs,
	\newblock Springer, New York, 
	\newblock (2007)	
	
\bibitem{Kempe}
	J.\,Kempe,
	\newblock  Quantum random walks - an introductory overview,
	\newblock Contemporary Physics, {\bf 44}, pp. 307--327,
	\newblock (2003)
	
\bibitem{Kendon}
	V.\,M.\,Kendon, C.\,Tamon,
	\newblock Perfect state transfer in quantum walks on graphs,
	\newblock Journal of Computational and Theoretical Nanoscience {\bf 8}, pp. 422--433,
	\newblock (2011)	
	
\bibitem{Sato11}
	N.\,Konno, I.\,Sato,
	\newblock On the relation between quantum walks and zeta functions,
	\newblock Quantum Information Processing, {\bf 11}, pp. 341--349,
	\newblock (2011)  
	
\bibitem{KOS}
	N.\,Konno, N.\,Obata, E.\,Segawa,
	\newblock Localization of the Grover walks on spidernets and free Meixner laws,
	\newblock Communications in Mathematical Physics, {\bf 322}, pp. 667--695,	
	\newblock (2013)
	
\bibitem{Konno14}
	N.\,Konno,
	\newblock Quantum Walk,
	\newblock Morikita Publishing Corporation Itd,
	\newblock (2014)	
	
\bibitem{Konno13}
  	N.\,Konno, Y.\,Shimizu, M.\,Takei, 
  	\newblock Periodicity for the Hadamard walk on cycle,
  	\newblock Interdisciplinary Information Sciences, {\bf 23}, pp. 1--8,
  	\newblock (2017)
	
\bibitem{Wang}
	K.\,Manouchehri, J.\,Wang,
	\newblock Physical Implementation of Quantum Walks,
	\newblock Springer, Berlin,
	\newblock (2014)	
	
\bibitem{tri}
	B.\,N.\,Parlett,
	\newblock The Symmetric Eigenvalue Problem,
	\newblock Prentice-Hall, Inc. Upper Saddle River, NJ, USA, 
	\newblock (1998)
		
\bibitem{Portugal}
	R.\,Porugal,
	\newblock Quantum Walks and Search Algorithms,
	\newblock Springer, Berlin,
	\newblock (2013)  	 
	
\bibitem{Rojo05}
	O.\,Rojo, R.\,Soto,
	\newblock The spectra of the adjacency matrix and Laplacian matrix for some balanced trees,
	\newblock Liner algebra and its applications, {\bf 43}, pp. 97--117,
	\newblock (2005)
	
\bibitem{Rojo07}
	O.\,Rojo, M.\,Robbiano,
	\newblock  An explicit formula for eigenvalues of generalized Bethe trees and upper bounds on the largest eigenvalue of any tree,
	\newblock   Liner algebra and its applications, {\bf 427}, pp. 138--150,
	\newblock (2007)
	
\bibitem{Stefanak}
	M.\,Stefanak, S.\,Skoupy,
	\newblock Perfect state transfer by means of discrete-time quantum walk search algorithms on highly symmetric graphs,
	\newblock Physical Review A, {\bf 94}, 022301,  
	\newblock (2016) 	
	
\bibitem{Stra07}
	F.\,W.\,Strauch,
	\newblock Relativistic effects and rigorous limits for discrete and continuous-time quantum walks,
	\newblock Journal of Mathematical Physics, {\bf 48}, 082102, 
	\newblock (2007)		
	
\bibitem{Sze04}
	 M.\,Szegedy,
	 \newblock Quantum speed-up of Marcov chain based algorithms,
     \newblock Proceedings of the 45th IEEE Symposium on Foundations of Computer Science, pp. 32--41,
	 \newblock (2004)  	
	 
\bibitem{Yoshie}
	Y.\,Yoshie,
	\newblock A characterization of the graphs to induce periodic Grover walk,
	\newblock accepted to publication to Yokohama Mathematical Journal,
	\newblock (2017)	 
  
  \end{thebibliography}
\end{document}